# MÉTODO DE AVALIAÇÃO APROXIMADA PARA A PROBABILIDADE DA UNIÃO DE EVENTOS INDEPENDENTES

# APPROXIMATE EVALUATION METHOD FOR THE PROBABILITY OF THE UNION OF INDEPENDENT EVENTS

EDSON LUIZ URSINI, PAULO S. MARTINS
FACULDADE DE TECNOLOGIA – UNICAMP
*{Ursini,* Paulo}@ft.unicamp.br

*Resumo – O cálculo da probabilidade de união de um grande número de eventos independentes requer várias combinações envolvendo o cálculo fatorial e exigindo o uso de computadores de alto desempenho e várias horas de processamento. Os limites da probabilidade de união e simplificações no seu cálculo são úteis na análise de problemas estocásticos em várias áreas como confiabilidade dos sistemas, sistemas biológicos, sistemas de tempo real com tolerância a falhas, teoria da probabilidade, teoria da informação, e comunicações dentre outros. Propomos uma aproximação que obtém a probabilidade da união de vários eventos independentes utilizando a média aritmética da probabilidade de todos eles. Os resultados aproximados têm um erro muito próximo, mas maior do que o erro verdadeiro. Isso permite um número muito menor de operações com um resultado semelhante e com maior simplicidade.*

*Palavras-chave: Probabilidade de União. Limites da Probabilidade de União. Eventos Independentes. Avaliação Aproximada.*

*Abstract – The evaluation of the probability of union of a large number of independent events requires several combinations involving the factorial and the use of high performance computers with several hours of processing. Bounds and simplifications on the probability of the union are useful in the analysis of stochastic problems across various areas including (but not limited to) systems reliability, biological systems, real-time fault-tolerant systems, probability theory, information theory and communications. We propose an approximation to evaluate the probability of the union of several independent events that uses the arithmetic mean of the probability of all of them. The approximate results are very close to, but larger than the exact values. The method allows a much smaller number of operations with a similar result and more simplicity.*

*Keywords: Probability of Union. Bounds on the Probability, Independent Events. Approximate Evaluation.*

## I. INTRODUCTION

Most systems are structured around a number of components, elements, parts or devices which as a whole provide the overall behavior. In this work, we use the word "system" in a holistic sense to encompass a large class of physical and non-physical systems such as mechanical (e.g. automobiles and aircraft), biological, computational (software and hardware), information systems, chains of command in military organizations, cyber-physical and real-time embedded systems among others. If we take one arbitrary but nevertheless relevant property of such systems, for example their reliability (for the sake of argument), its overall system reliability may be obtained from the individual reliability of its components. In particular, if the system components are arranged in series, failure of one component may compromise the overall system's reliability (clearly, if no form of redundancy is provided).

Within this context, the probability of failure of a system consisting of *n* components in series is given by the probability of the union of the event of each component failing independently. However, the calculus of the probability of the union of a large number of independent and non-mutually exclusive events is computationally intensive, as it requires several combinations each involving the factorial. For example, in the case of 300 devices that may fail we need to resort to combinations of 300 taken *k* at a time, where *k* may vary from 1 to 300. Therefore, the goal of our work is to provide a method that simplifies this calculus by avoiding the calculation of all of these combinations, i.e. the proposed method relies on an approximate calculus while also evaluating the error incurred by means of bounds of the probability of the union. More specifically, *in this work we aim at calculating the approximate value of the probability of the union of independent events by using the mean value of the probability of the events*. To our knowledge, we have found no work in the literature that has adopted this approach to the estimation of the probability of the union.

Notice that, although the main focus here is to show in a didactic way the application of the method to devices or components with estimated fault probabilities (e. g. resistors or devices in series such as in the realm of reliability theory or several serial links in a telecommunication system), our work accommodates far more general scenarios of the union of independent events than the one covered here.

One of the earliest papers on the topic is the work by Miller (1968), who focuses on issues related to processing time and memory space. However, unlike our work, in the author's traditional approach there is no reference to the use of the mean value to calculate the probability of the union of events. Legg et al (2011) address the problem of calculating a very large number of independent events, but they use true (real) values (i.e. not the average value as we propose in this work) which implies a very long processing time. For up to 100 items, their convergence is very fast. However, time grows exponentially beyond 100 devices. As an example, the calculation performed with 1000 devices (items) took more than 17 hours. Very often the processing time is a stringent requirement and a shorter processing time may be more important than exact precision. The last two papers specifically addressed the probability of the union of independent events.

The paper from Caen (1997) presents a lower bound on the probability of the union using only the individual probabilities and the probabilities of the joint events, two by two. In another work, Prékopa and Gao (2005) use two linear programming boundaries to establish another boundary for the probability of union of events since there is no complete knowledge of the joint probability of all events. The authors also generalize the work by De Caen (1997). The paper from Kuai et al (2000) also seeks to use the individual probabilities and the probability of the union of events, two by two, to establish the bounds for the probability of the union. The work by Veneziani (2009) improves upon the work of Prékopa and Gao (2005) by the inclusion of weights for the resolution of the linear system. Our approach is simpler because it assumes the independence between events. Furthermore, it seeks an approximate calculation by using the mean of the probability of occurrence of individual events.

The work of Kounias (1995) uses the Bonferroni inequality or the Poisson approximation to evaluate the error resulting from not knowing the joint probability of the events that are not independent. On the other hand, the paper by Yang et al (2016) derives lower bounds on the finite probability of union in terms of the individual event probabilities and a weighted sum of the pairwise event probabilities that have at most pseudo-polynomial computational complexity. This work also generalizes some recent works. Hollenback and Moss (2011) address the issue of approximation, and offer the calculus of the limits to the probabilities. However, unlike our work, they do not deal with the issue of using the mean value of probabilities. The objective of our work is mainly to attempt to cut down computational overheads such as processor time and memory. Concerning the independence of events, our approach may also be useful in cases where the events are approximately independent.

The remainder of this paper is organized as follows: Section II presents the background on the calculus of the probability of independent events; Section III specializes this calculus for the case where the devices have the same probability of failure; Section IV derives the incurred error when the devices have different probability of failure; Section V deals with the issue of the incurred error when a smaller number of terms are adopted; In Section VI we propose the use of the mean value of the probability of failure and we evaluate its corresponding approximation errors. Finally, Section VII addresses our remarks and conclusions.

## II. PROBABILITY OF THE UNION

We consider a system with $n$ components (or devices) in series, for didactic purpose, $n$ resistors or $n$ telecommunication links (the same could be extended for probability of union of any independent sets). The events $A_i, i = 1, \cdots, n$ correspond to the failure operation of each of the components. Thus, the failure probability of each component is $P(A_i) = p_i$.

In turn, $P(\overline{A_i}) = 1 - P(A_i) = 1 - p_i$ is the probability that the device $i$ do not fail. The probability of the union for $n$ independent devices, each with a probability of failure $p_i$, is given by (GAVIN, 2016):

$$P(A_1 \cup A_2 \cup \cdots \cup A_n) = 1 - P(\overline{A_1 \cap A_2 \cap \cdots \cap A_n}) =$$
$$= 1 - P(\overline{A_1})P(\overline{A_2})\cdots P(\overline{A_n}) = \quad (1)$$
$$= 1 - (1 - p_1)(1 - p_2)\cdots(1 - p_n)$$

In the preceding expansion the De Morgan theorem was employed along with the knowledge that the events are independent (i.e. devices fail or survive independently of one another).

## III. PROBABILITY OF THE UNION FOR EVENTS WITH THE SAME PROBABILITY

If the events have the same probability, then:

$$P(A_1 \cup A_2 \cup \ldots \cup A_n) = 1 - (1 - p)^n \quad (2)$$

Equation (2) alone is not sufficient to find out the number of operations required to calculate the final probability. Thus, by further developing expression (2) for equal probabilities we have:

$$P(A_1 \cup A_2 \cup \cdots \cup A_n) = \binom{n}{1}p - \binom{n}{2}p^2 + \binom{n}{3}p^3 - \cdots (-1)^n \binom{n}{n}p^n$$

.

The terms alternate their signal, starting from a positive sign. The first term corresponds to the case where the probabilities of failure are mutually exclusive. The positive sign (+) in the last term occurs for an odd $n$ and the negative (-) for an even value of $n$. The right-hand side of this summation may be written as:

$$\sum_{i=1}^{n}(-1)^{i-1}\binom{n}{i}p^i \quad (3)$$

However, we may expand this expression:

$$\sum_{i=1}^{n}(-1)^{i-1}\binom{n}{i}p^i = \sum_{i=1}^{n}-\binom{n}{i}(-p)^i = -\sum_{i=1}^{n}\binom{n}{i}(-p)^i \quad (4)$$

The summation $\sum_{i=0}^{n}\binom{n}{i}(-p)^i$ is the Newton's binomial which is equal to $(1-p)^n$. Thus, observing that in (3) we are missing only the first term of the binomial (which is 1),

$$-\sum_{i=1}^{n}\binom{n}{i}(-p)^i = -[(1-p)^n - 1] = 1 - (1-p)^n$$

Therefore, the probability of the union for $n$ devices is obtained by different means as shown in (2):

$$P(A_1 \cup A_2 \cup \ldots \cup A_n) = 1 - (1-p)^n$$

This is the correct value for finding the probability of the union for $n$ statistically independent devices with the same probability of failure. On the other hand, if we want to use fewer terms, we must truncate equation (4) by

$$-\sum_{i=1}^{m}\binom{n}{i}(-p)^{i} \quad m < n \quad (4a).$$

If we carried out the summation in $n$, we would obtain equation (2)

## IV. PROBABILITY OF THE UNION FOR EVENTS WITH DIFFERENT PROBABILITIES

When the probabilities of the components differ from each other, equation (3) has to be recast as:

$$P(A_1 \cup A_2 \cup \cdots \cup A_n) = \\ (p_1 + \cdots + p_n) - (p_1 p_2 + p_1 p_3 + \cdots + p_{n-1} p_n) + \\ (p_1 p_2 p_3 + p_1 p_2 p_4 + \cdots + p_{n-2} p_{n-1} p_n) - \cdots \\ + (-1)^n (p_1 \cdots p_{n-2} p_{n-1} p_n) \quad (5)$$

Each term corresponding to each parenthesis in (5) has $\binom{n}{i}$ combinations. The first term has $\binom{n}{1} = n$ values and it corresponds to the particular case where the probabilities of failure are mutually exclusive. The second term has $\binom{n}{2} = \frac{n(n-1)}{2!}$ products of two probability values. The third term has $\binom{n}{3} = \frac{n(n-1)(n-2)}{3!}$ products of three probability values, and so on, until $\binom{n}{n} = 1$ one product of $n$ probability values. All terms should be multiplied by $(-1)^{i-1}$.

## V. INCURRED ERROR

In order to illustrate the concept, we now present the complete calculus for two and three independent devices. For 2 devices we have: $(p_1 + p_2) - (p_1 p_2)$ or $2p - p^2$ for equal probabilities; for 3 devices we have:

$$(p_1 + p_2 + p_3) - (p_1 p_2 + p_1 p_3 + p_2 p_3) + (p_1 p_2 p_3) \quad \text{or}$$
$3p - 3p^2 + p^3$ when the probabilities are the same.

For a large number of devices the calculus of probability from equation (4) is a challenge due to the large number of combinations. Furthermore, the complete numeric calculus of the probability of the union of all devices requires substantial processing power. This task may be simplified by employing (3). If it is possible to find the probability of failure from (3) and the number of terms, for an error set a priori from (1), the resulting error is controlled by the user (in addition to allow efficient algorithms with relatively small approximations).

In order to evaluate the error incurred by neglecting the subsequent terms (in percentage), we apply the modulus of the percentage value obtained in a truncating point $m < n$:

$$\left| \frac{P(A_1 \cup A_2 \cup \ldots \cup A_n) - \sum_{i=1}^{m}\binom{n}{i}\text{SUM PROB}(i)}{P(A_1 \cup A_2 \cup \ldots \cup A_n)} \right| \quad m < n \quad (6)$$

SUM PROB($i$) in equation (5) means the following: $i = 1$ corresponds to the particular case where the probabilities of failure are mutually exclusive, if $i = 2$, there is a product of two probability values, and if $i = 3$ there is a product of three probability values and so on until $m$ ($m < n$). If $m = n$, the calculation is complete and equation (5) is the same as equation (1) (if the probabilities are equal, SUM PROB($i$) is the same as in equation (4a).

Notice that the computation of all the terms in the equation above ($m = n$) leads to a null/zero error. A numeric example, using equation (1), is the case for two devices with $p_1 = 0.1$ and $p_2 = 0.3$ where

$$P(A_1 \cup A_2) = 1 - (1 - p_1)(1 - p_2) = p_1 + p_2 - p_1 p_2 = \\ = 0.1 + 0.3 - 0.03 = 0.37$$

## VI. PROPOSED APPROACH

Assuming devices with the same probability of failure, the mean probability $\overline{p}$ of failure for each device is $\overline{p} = \frac{p_1 + p_2}{2} = 0.2$. Thus, the probability of failure of the set for the mean value $\overline{p} = 0.2$, and using equation (2), is

$$P(A_1 \cup A_2) = 1 - (1 - 0.2)^2 = 0.2 + 0.2 - 0.04 = 0.36.$$

The approximation (truncation) in this case can only be applied to the first term (it falls under the special case where we assume that the events are mutually exclusive). This is a special case, but shows an approximate error of 11% (due to using a mean value) and a correct error of 8% when using only one term.

Table 1 introduces five examples of complete and approximate probability of failure (using the mean of failures and with all the terms).

Table 1 - Approximate and exact values of the probability of the union (failure).

| # of Devices n | Probability of failure | Probability (mean) ($\overline{p}$) | Exact value equation (1) | Approxi-mate value equation (2) |
|---|---|---|---|---|
| 2 | 0.1 and 0.3 | 0.20 | 0.3700 | 0.3600 |
| 3 | 0.1, 0.3 and 0.5 | 0.30 | 0.6850 | 0.6570 |
| 4 | 0.1, 0.2, 0.2 and 0.3 | 0.20 | 0.5968 | 0.5904 |
| 4 | 0.5, 0.8, 0.2 and 0.4 | 0.48 | 0.9520 | 0.9240 |
| 5 | 0.1, 0.2, 0.2, 0.3 and 0.2 | 0.20 | 0.6774 | 0.6723 |

Table 2 illustrates the resulting errors for *n*-1 terms (term 1 corresponds to a simple probability, term 2 are probabilities represented as two-by-two products, ..., and so forth) for finding the probability of the union. Assuming *n* devices with known probability of failure, the exact error that results when using *n* - 1 terms is given by equations (5) and (6) (with 2 devices we have *n* - 1 = 1, a single term; for three devices we have two terms, and so on). For the approximate error we use equations (2) and (6) because we work with the mean value, $\bar{p}$, so the error that would incur for the mean value is equation (2) for $1-(1-\bar{p})^n$:

Table 2 - Incurred error values (%) for *n*-1 terms in equation (5) or in equation (4).

| # of devices *n* | Probability of failure | Probability (mean) ($\bar{p}$) | Exact error (%) eq.(5 & 6) | Approx. error (%) eq. (4 &6) |
|---|---|---|---|---|
| 2 | 0.1 and 0.3 | 0.20 | one term 8.0 | one term 11.0 |
| 3 | 0.1, 0.3 and 0.5 | 0.30 | two terms 2.0 | two terms 4.0 |
| 4 | 0.1, 0.2, 0.2 and 0.3 | 0.20 | three terms 0.2 | three terms 0.27 |
| 4 | 0.5, 0.8, 0.2 and 0.4 | 0.48 | three terms 3.36 | three terms 5.5 |
| 5 | 0.1, 0.2, 0.2, 0.3 and 0.2 | 0.20 | four terms $3.5 \cdot 10^{-4}$ | four terms $4.8 \cdot 10^{-4}$ |

Remarks:

(i) The errors of the approximated probabilities are very close to the exact ones, although always larger. This ensures that the actual error when connecting the actual operations will be less than the value obtained when calculating with the average probability because the operations obey the actual values (and not the mean values);
(ii) In general, the examples for the probability values 0.1, 0.2, etc., are already relatively large values for practical, real-world scenarios dealing with failures;
(iii) The larger the number of items (i.e. devices or terms in the approximation) the smaller the incurred error.

*6.1 - Incurred error as a function of the number of terms considered*

We analyze the evolution of the incurred errors, the difference between the exact values (i.e. the "correct value") in comparison to the approximate values. In any case, the approximate value is close to the real value but always worse than the former. Using the five examples from Table I, we show the evolution of the errors incurred term by term. For two components, Fig. 1 shows the estimated error as a function of the number of terms. Notice that the *x* and *y*-axis have continuous values only for illustration purposes. However, these values are discrete and they correspond to the number of terms considered in the approximation. The Appendix A shows the Matlab code for calculating the approximate value (or the exact value if the items have the same probability of failure). The function ([error] = Prob_union_v2 (p, m, n]) has as input the values of *p*, *m* (truncation) and *n*, and it returns the incurred error.

Figure 1- Estimated error for two devices

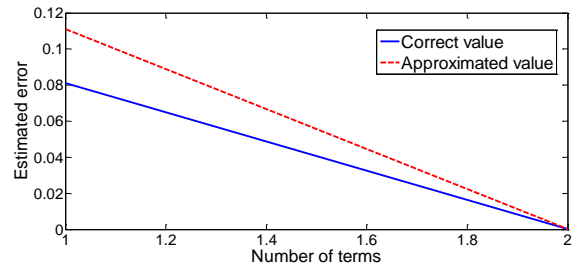

Consider an example for 3 devices with $p_1 = 0.1$, $p_2 = 0.3$ and $p_3 = 0.5$. Assuming that the devices have the same probability of failure, the mean value is $\bar{p} = \frac{p_1 + p_2 + p_3}{3} = 0.3$. Fig. 2 shows that the correct error for a single term is 32% and the approximate error is 37%. With two terms, the correct value is 2% and the approximate value is 4%.

Figure 2 - Estimated error for three devices

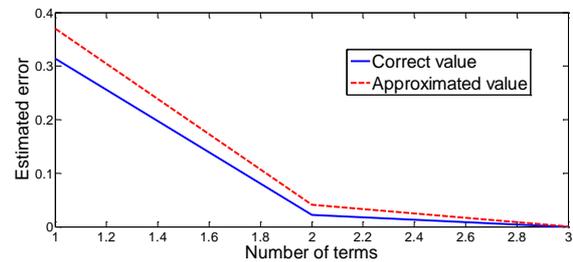

The following example utilizes 4 devices and considers that $p_1 = 0.1$, $p_2 = 0.2$, $p_3 = 0.2$ and $p_4 = 0.3$ (fourth row in Tables I and II). Assuming devices with the same probability of failure, then $\bar{p} = \frac{p_1 + p_2 + p_3 + p_4}{4} = 0.2$. Fig. 3 shows that for one term the correct error is 34% and the approximated error is 36%. For two terms the correct error is 4.5% and the approximated value is 5%.

Figure 3 - Estimation error for four devices

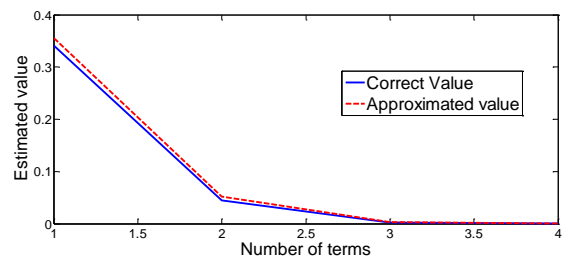

In Tables I and II another example is included considering $p_1 = 0.5$, $p_2 = 0.8$, $p_3 = 0.2$ and $p_4 = 0.4$. This example was just to illustrate the fact that the error grows as we increase the difference between the probabilities. Nevertheless, these values are not common in component failure and their difference is not common either.

The next example employs 5 devices considering $p_1 = 0.1$, $p_2 = 0.2$, $p_3 = 0.2$, $p_4 = 0.3$ and $p_5 = 0.2$. Assuming devices with the same probability of failure $p_i$ we have $\overline{p} = \dfrac{p_1 + p_2 + p_3 + p_4 + p_5}{5} = 0.2$. Fig. 4 shows that for just one term the correct error is 47% and the approximated one is 48%. With only two terms the exact error is 10% and the approximate one is 11%. For three terms both the exact and the approximate error values are only 1%.

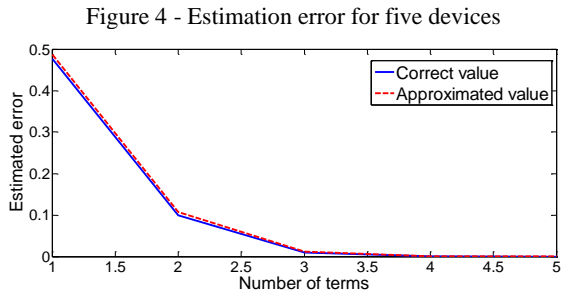

Figure 4 - Estimation error for five devices

An important fact is that the approximate values are always "worse" than the real values (although close to). This ensures that the correct result will be better than (because the values, in fact, will be the real ones and not the average ones) stopping condition imposed on the algorithm that implements the proposed method (i.e. using the approximated calculus furnished by equations (4) and (6)). Table 3 illustrates an example using 100 items subject to failure for two different values of mean probability of failure (0.1 and 0.01). The values of the maximum incurred error as a function of the number of terms in equation (4) may be evaluated.

Table 3 - Max. error as a function of the number of terms for 100 devices

| $\overline{p}$ mean | # of terms | Maximum Error (%) |
|---|---|---|
| 0.1 | 26 | 0.15 |
| 0.1 | 27 | 0.04 |
| 0.1 | 28 | 0.01 |
| 0.1 | 29 | 0.002 |
| 0.1 | 30 | 0.0005 |
| 0.01 | 1 | 57.74 |
| 0.01 | 2 | 20.34 |
| 0.01 | 3 | 5.16 |
| 0.01 | 4 | 1.02 |
| 0.01 | 5 | 0.17 |
| 0.01 | 6 | 0.023 |

An observation about the calculations is that the maximum value of the factorial in the processor in which they were made is $n! = 170$. For values larger than 170 (for instance to calculate up to 300 devices), we may apply the Stirling's formula. However, in this case it is important that we perform a deeper analysis of the incurred error, which is deferred to a future work. When the probabilities involved are very large, using a few terms can lead to erroneous results, of the type with probability greater than one. In this case, it is important to use a program ([error, i] = Prob_union_rev2(p,n,re) in Appendix B) that calculates the number of terms (i) for the minimum required error (re).

The Matlab algorithm which implements the method proposed in this work, with respect to the maximum permissible error, is found in the Appendix B.

## VII. REMARKS AND CONCLUSIONS

The estimation of the overall probability of union of events (such as in the reliability of telecommunication systems) relies on the calculus of the probability of the union of events - and it may require a high performance computational power. This work provided an alternative method that simplifies the calculation of the system probability of the union (of failure) concerning independent events, as it relies on a single value, i.e. the mean probability of failure for the devices, instead of the more traditional equations which often lead to a combinatorial explosion (and therefore high computational cost) particularly while estimating the probability of the union (failure) in large systems.

The approximate method allows that algorithms that correctly calculate the probability of $n$ independent events have a stop criterion that ensures a controlled error. In addition to the fact that the approximate calculus is much simpler to execute (equation (2) and mainly equation (4)), which is in itself a significant gain, it allows the calculation of probabilities of a large number of devices with a high precision. The approximate evaluation leads to probabilities that are larger than the values from the complete calculation, and this conservative behavior is appropriate in many areas such as safety, reliability and real-time systems. The actual values are smaller than these values from the approximation because the devices are those with real probabilities and not those with average probabilities. We must also consider that, if the probabilities of failure are too close (i.e. if devices have the same probability) the result is still better. However, it is always possible to use the terms of equation (4) instead of the ones from equation (5), i.e. to calculate the error using equation (4) knowing that the true error is smaller since it follows equation (5). The method may be used as a planning tool and to quickly (also with more simplicity) dimension the probability of failure.

The simplicity of the algorithm may be contemplated in the Appendix, where we present the Matlab code that implements the proposed calculus. Another fact worth pointing out is that the approximation is also strengthened by the fact that the calculus of the probability of system failure (or probability of the union of independent events) is intrinsically imprecise due to the estimation error of the individual (i.e. component) probabilities.

As future work, the approach could be extended to minimize operations in other types of configurations, e.g. parallel or serial-parallel devices, always aiming to use less computational resources such as processing time or memory allocation. As a last consideration, the approach could also be employed in the case where the probabilities of the involved components are nearly independent.

## IX. COPYRIGHT

APPENDIX:

A) Matlab algorithm depending on the number of terms

```
function [error] = Prob_union_v2(p,m,n)

% m is the number of terms considered
% n is the number of devices
% error is the numerical calculation mistake
s=0;
for i=1:m,
s=s+(-1)^(i-1)*prod([1:n])/(prod([1:i])*prod([1:n-i]))*p^i;% part of equation (4a)
end
prob_v = 1-(1-p)^n; % equation (2)
error = abs(s-prob_v)/prob_v; % equation (6)
```

B) Matlab algorithm depending on required error

```
function [error,i] = Prob_union_rev2(p,n,re)

% n is the number of devices
% p is the specified probability
% re is the specified error
% error is the returned error value
% i is the number of terms considered
s=0; error=1;
i=1;
while  error > re;
s=s+(-1)^(i-1)*prod([1:n])/(prod([1:i])*prod([1:n-i]))*p^i;% part of equation (4a)
prob_v=1-(1-p)^n;% equation (2)
error=abs(s-prob_v)/prob_v;% equation (6)
i=i+1;
end
i=i-1
```